\documentclass[11pt]{article}
\usepackage[T1]{fontenc}
\usepackage[latin9]{inputenc}
\usepackage{amsmath}
\usepackage{amssymb}
\usepackage{esint}

\makeatletter
\usepackage{amsfonts}%TCIDATA{OutputFilter=latex2.dll}
\newtheorem{theorem}{Theorem}

\makeatother

\begin{document}
\begin{center}
{\Large Regularity Criteria of BKM type in Distributional Spaces for
the 3-D Navier-Stokes Equations on Bounded Domains}
\par\end{center}{\Large \par}

\begin{center}
{\Large \vspace{0.2in}
}
\par\end{center}{\Large \par}

\begin{center}
{\Large Joel Avrin }
\par\end{center}{\Large \par}

\begin{center}
{\Large Department of Mathematics and Statistics}
\par\end{center}{\Large \par}

\begin{center}
{\Large University of North Carolina at Charlotte}
\par\end{center}{\Large \par}

\begin{center}
{\Large jdavrin@uncc.edu\vspace{0.2in}
}
\par\end{center}{\Large \par}

\begin{center}
\underline{Abstract} 
\par\end{center}

In the classic work of Beale-Kato-Majda ({[}2{]}) for the Euler equations
in $\mathbb{R^{\mathrm{3}}}$, regularity of a solution throughout
a given interval $[0,T_{*}]$ is obtained provided that the curl $\omega$
satisfies $\omega\in L^{1}((0,T);L^{\infty}(\mathbb{R^{\textrm{\ensuremath{3}}}})$
for all $T<T_{*}$, and the authors noted that the arguments apply
equally well to the Navier-Stokes equations (NSE) in $\mathbb{R^{\mathrm{3}}}$.
In later works by various authors the spatial $L^{\infty}$-criterion
imposed on the curl was generalized to a $BMO$ criterion, and later
to a Besov space criterion, in both the Euler and NSE cases ({[}9{]},
{[}10{]}, {[}11{]}). Meanwhile, the authors in {[}2{]} remarked that
additional ideas seem necessary to obtain results of this type on
bounded spatial domains. Efforts in this direction in {[}8{]} for
the NSE case produced regularity results with the $BMO$ criterion
imposed on localized balls. 

In this paper for the NSE case and on general bounded domains $\Omega$
in $\mathbb{R^{\mathrm{3}}}$, we obtain a regularity result of BKM
type that goes beyond function spaces to spatially allow $\omega$
to be a distribution. This is done by making a new connection between
a well-known vector calculus result and the clasical regularity criteria
of Serrin type ({[}12{]}, {[}14{]}, {[}15{]}, {[}18{]}). Specifically,
for certain Sobolev spaces $H^{s,p}(\Omega)$ suitably defined for
$s<0$ we show that if $u$ is a Leray solution of the 3-D NSE on
the interval $(0,T)$ and if $\omega\in L^{s}((0,T);H^{-1,p}(\Omega))$
where $\frac{2}{s}+\frac{3}{p}=1$ for some $p\in(3,\infty]$, then
$u$ is a regular solution on $\left(0,T]\right)$; in particular
for $p=\infty$ we have a regular solution when $\omega\in L^{2}((0,T);H^{-1,\infty}(\Omega))$,
which directly strengthens the results in {[}2{]} by one order of
(negative) derivative in terms of the spatial criteria for regularity.
Our results thus impose more stringent conditions on time than the
BKM results and their generalizations described above, but as far
as we are aware the results here represent the first of BKM type for
the NSE that allow $\omega$ to spatially be a distribution.


\section{Introduction}

We consider the 3-D Navier-Stokes equations for viscous incompressible
homogeneous flow
\begin{align}
u{}_{t}+\nu Au+\left(u\cdot\nabla\right)u+\nabla p=g,\tag{1.1a}\nonumber \\
\nabla\bullet u=0.\tag{1.1b}
\end{align}
Here $\Omega$ is a bounded spatial domain in $\mathbb{R}^{3}$ with
sufficiently smooth boundary and $u=\left(u_{1},u_{2},u_{3}\right)$
with $u_{i}=u_{i}\left(x,t\right),x\in\Omega,1\leq i\leq n$ and $t\geq0$.
The external force is $g=\left(g_{1},g_{2},g_{3}\right)$, with $g_{i}=g_{i}\left(x,t\right)$,
and $p=p\left(x,t\right)$ is the pressure. The domain $\Omega$ can
be either a periodic box or a Lipschitz domain with zero (no-slip)
boundary conditions; in the latter case, or by ''moding out'' the
constant vectors as in standard practice in the former case, $A=-\Delta$
has eigenvalues $0<\lambda_{1}<\lambda_{2}<\cdots$ with corresponding
eigenspaces $E_{1},E_{2},\cdots$, so that in particular $A$ is a
positive definite operator and $A^{-1}$ is a well-defined bounded
operator on the Banach spaces $L^{p}(\Omega)$, $p\in[1,\infty)$.
Let $\partial_{x}$ denote the operator $\frac{\partial}{\partial x}$
then with similar definitions for $\partial_{y}$ and $\partial_{z}$
we have that $\nabla\bullet u=div\text{ }u=\partial_{x}u_{1}+\partial_{y}u_{2}+\partial_{z}u_{3}$.
Of particular interest also is the curl $\omega$ defined by $\omega=\nabla\times u=(\partial_{y}u_{3}-\partial_{z}u_{2},\partial_{z}u_{1}-\partial_{x}u_{3},\partial_{x}u_{2}-\partial_{y}u_{1})$.
With zero viscosity ($\nu=0$) the system (1.1) becomes the Euler
system
\begin{align}
u_{t}+\left(u\cdot\nabla\right)u+\nabla p=0,\tag{1.2a}\nonumber \\
\nabla\bullet u=0.\tag{1.2b}
\end{align}
In the classical work of Beale/Kato/Majda ({[}2{]}), $\Omega=\mathbb{\mathbb{R^{\textrm{\ensuremath{3}}}}}$,
$g=0$, and regularity for a smooth solution of (1.2) throughout a
given interval $[0,T_{*}]$ is obtained provided that $\omega\in L^{1}((0,T);L^{\infty}(\mathbb{R^{\textrm{\ensuremath{3}}}})$
for all $T<T_{*}$. Central to the arguments in {[}2{]} is the formula
$u=-\nabla\times(\nabla^{-1}\omega)$ which in $\mathbb{R^{\textrm{\ensuremath{3}}}}$
is given explicitly by appropriately available kernels via the Biot-Savart
Law. The authors note that the results hold for periodic flow with
minor modification, and they note that the results apply to the NSE
as well. 

The results in {[}2{]} for (1.1) were extended in the case $\Omega=\mathbb{R^{\textrm{\ensuremath{n}}}}$
in {[}9{]} to allow $\omega\in L^{2}((0,T);BMO)$ where $BMO$ denotes
the class of functions of bounded mean oscillation. Later in {[}10{]}
this condition was extended to both (1.1) and (1.2) to allow $\omega\in L^{1}((0,T);BMO)$.
Then in {[}11{]} the results in {[}10{]} were extended to allow $BMO$
to be replaced by the Besov space $B_{\infty,\infty}^{0}$. The same
regularity criterion developed in {[}11{]} was then established in
the case $n=3$ for the Boussinesq system, the MHD system, and a fluid
system with the linear Soret effect in {[}4{]}, {[}13{]}, and {[}5{]},
respectively. Meanwhile the authors of {[}2{]} noted that a more involved
proof using additional ideas seems necessary for bounded spatial domains.
In {[}8{]} regularity results were obtained for the NSE case by imposing
the $BMO$ condition on localized balls. 

In this paper for the NSE case and on general bounded domains $\Omega$
in $\mathbb{R^{\textrm{\ensuremath{3}}}}$ with sufficiently smooth
boundary we will obtain regularity results of BKM type in which the
spatial criteria that we impose on $\omega$ will allow $\omega$
to lie in negative Sobolev spaces. Thus a.e. for each $t$ the curl
$\omega(\bullet,t)$ is allowed to be a distribution. 

In proving our results we will make use of the classic regularity
criteria for the Navier-Stokes equations which establish regularity
of Leray solutions (see the definition in section 2 below) provided
that $u\in L^{\theta}((0,T);L^{p}(\Omega)$ and $\theta,p,n$ satisfy
$\frac{2}{\theta}+\frac{n}{p}=2$, $n<p\leq\infty$. Here $\Omega=\mathbb{R^{\mathrm{\textrm{n}}}}$
or under suitable conditions such as those assumed here $\Omega$
is a bounded domain; see {[}12{]}, {[}14{]}, {[}15{]}, {[}18{]}, and
the references contained therein. Preliminary results toward extending
these classic results to the borderline case $n=p$ were obtained
in {[}6{]}, {[}7{]}, {[}19{]}, {[}20{]} (see also the references contained
therein), and recently this borderline result was obtained in the
case $\Omega=\mathbb{\mathcal{\mathbb{R^{\textrm{3}}}}}$ ({[}3{]},
{[}16{]}). It is as yet unknown if the borderline case can be obtained
on bounded domains. 

The other main component used in establishing our results will be
the well-known vector-calculus identity 
\begin{equation}
Av=\nabla\times\nabla\times v\tag{1.3}
\end{equation}
holding for smooth divergence-free vector fields on $\Omega$. The
smoothness we require for the boundary of $\Omega$ is that the usual
Sobolev inequalities hold. Since on $\Omega$ under these conditions
and for the assumed boundary conditions (e.g. zero Dirichlet) the
operator $A$ is invertible, we have from (1.3) that $v=A^{-1}(\nabla\times\nabla\times v)$
which provides an alternative relationship between $u$ and $\omega$
similar to $u=-\nabla\times(\nabla^{-1}\omega)$ but more adaptable
to bounded domains and more directly applicable to our techniques.
The following result easily generalizes the identity $v=A^{-1}(\nabla\times\nabla\times v)$: 

\begin{theorem} If $u$ is a smooth enough solution of (1.1) or (1.2)
then $A^{s}u=A^{s-1}(\nabla\times\omega)$, where $s$ is any order
allowed by the smoothness of $u$. \end{theorem}

Here as in standard fashion we let $H_{0}\equiv\overline{\{v\in C_{0}^{\infty}(\Omega):\nabla\bullet v=0\}}_{L^{2}(\Omega)}$,
ie. $H_{0}$ is the closure in $L^{2}(\Omega)$ of the smooth compactly-supported
solenoidal vector fields. Note that $v$ need not be smooth in order
for the relationship $v=A^{-1}(\nabla\times\nabla\times v)$ to hold,
since $A^{-1}(\nabla\times\nabla\times\bullet)$ defines a bounded
operator on $H_{0}$ as can bwe quickly seen (see section 2 below). 

Recall that the standard Sobolev spaces $W^{k,p}(\Omega)$ are defined
as $W^{k,p}(\Omega)\equiv\{v\in L^{p}(\Omega)\mid D^{\alpha}v\in L^{p}(\Omega)\forall|\alpha|\leq k\}$.
The corresponding negative Sobolev spaces are defined for each $k$
as the dual spaces of $W^{k,p'}(\Omega)$ where $\frac{1}{p}+\frac{1}{p'}=1$,
i.e. $W^{-k,p}(\Omega)\equiv(W^{k,p'}(\Omega))'$. Our use of Sobolev
spaces of negative or positive order is motivated by the following
characterization (see e.g. {[}17{]}):

\begin{theorem} Let $v\in\mathcal{D}'(\Omega)$, then $v\in W^{-k,p}(\Omega)$
if and only if $v={\displaystyle {\displaystyle \sum_{|\alpha|\leq k}D^{\alpha}w_{\alpha}}}$
where $w_{\alpha}\in L^{P}(\Omega)$. \end{theorem}

Theorem 2 says in a sense that $v\in W^{-k,p}(\Omega)$ iff $D^{-k}v\in L^{p}(\Omega)$;
we can make this more precise by defining suitable Sobolev spaces
$H^{-s,p}(\Omega)$ for any $s\geq0$ by $H^{-s,p}(\Omega)\equiv\{v\in\mathcal{D}'(\Omega)\mid A^{-s/2}v\in L^{p}(\Omega)\}$
where as noted we have assumed that $\Omega$ has a boundary smooth
enough such that the usual Sobolev spaces as well as the operators
$A^{-s/2}$ are well-defined. In fact, if we take this together with
the definition $H^{s,p}(\Omega)\equiv\{v\in D(A^{s/2})\}$ with norm
$\Vert v\Vert_{s,p}\equiv\Vert A^{s/2}v\Vert_{p}$ then we have a
consistent definition of $H^{s,p}(\Omega)$ for any real $s$ and
any $p\in(1,\infty]$. These are the Sobolev spaces we will work with,
of both negative and positive order. Similar spaces were defined and
used in {[}7{]} and {[}20{]} wherein $A$ in those cases was the Stokes
operator $-P\triangle$ where $P$ is the Leray projection onto the
solenoidal vectors. From the basic tools developed in Theorem 1, the
regularity criteria $u\in L^{\theta}((0,T);L^{p}(\Omega),\frac{2}{\theta}+\frac{3}{p}=2,3<p\leq\infty$
as noted above for the case $n=3$, and our definition here of the
spaces $H^{-s,p}(\Omega)$ we will establish the main result of this
paper:

\begin{theorem} Let $u$ be a Leray solution of (1.1) on the interval
$(0,T)$ and suppose that $\omega\in L^{\theta}((0,T);H^{-1,p}(\Omega))$
where $\frac{2}{\theta}+\frac{3}{p}=1$ for some $p\in(3,\infty]$.
Then $u$ can be continued to a regular solution of (1.1) on $\left(0,T]\right)$.
\end{theorem}

We remark that as an immediate corollary of Theorem 3 (and as noted
similarly in {[}2{]}, {[}9{]}, {[}10{]}, {[}11{]}) we have that if
the maximal existence time $T^{*}$ is finite then $\limsup_{t\uparrow T^{*}}\Vert\omega(t)\Vert_{H^{-1,\infty}(\Omega)}=\infty$.
Theorem 3 overlaps with the main result of {[}2{]} and the results
in {[}4{]}, {[}5{]}, {[}10{]}, {[}11{]}, {[}13{]} in that the condition
on the integrability in time is more stringent while the spatial requirement
on $\omega$ is more general. Specific to the case $p=\infty$ the
results in {[}2{]} require that $\omega\in L^{1}((0,T);L^{\infty}(\Omega))$
whereas here the corresponding condition is that $\omega\in L^{2}((0,T);H^{-1,\infty}(\Omega))$;
this means in particular that $A^{-1/2}\omega(\bullet,t)\in L^{\infty}(\Omega)$
a.e. for each $t$ in contrast with the requirements in {[}2{]} which
imply that $\omega(\bullet,t)\in L^{\infty}(\Omega)$ a.e. for each
$t$. Theorem 3 will follow by connecting the results of Theorem 1
with the regularity criteria $u\in L^{\theta}((0,T);L^{p}(\Omega),\frac{2}{\theta}+\frac{3}{p}=2,3<p\leq\infty$
by using a few duality arguments similar to those employed in {[}7{]}
and in {[}1{]}; af ter some preliminary discussion Theorems 1 \& 3
will be proven in the next section. In section 3 we will make some
concluding remarks and observations.

\section{Preliminaries and Proof of Theorem 3}

By a Leray solution of (1.1) on $(0,T)$ we mean a vector $u\in L^{\infty}((0,T);L^{2}(\Omega))\bigcap L^{2}((0,T);H^{1}(\Omega))$
satisfying, for each $v\in L^{\infty}((0,T);L^{2}(\Omega))\bigcap L^{2}((0,T);H^{1}(\Omega))$,
the equation 
\begin{equation}
\left(u(t),v\right)+\nu\int_{t_{0}}^{t}\left(A^{1/2}u,A^{1/2}v\right)+((u\cdot\nabla)u,v)ds=\left(u(t_{0}),v\right)+\int_{t_{0}}^{t}\left(g,v\right)ds\tag{2.1}
\end{equation}
 for all intervals $(t_{0},t)$ contained in $(0,T)$. Since $((u\cdot\nabla)u,uv)=-((\nabla\bullet u)u,u)=0$
and by the standard use of Young's inequality on the term $\left(u(t_{0}),u\right)$
we have by setting $v=u$ in (2.1) that 
\[
\frac{1}{2}\left\Vert u(t)\right\Vert _{2}^{2}+\nu\int_{t_{0}}^{t}\left\Vert A^{1/2}u\right\Vert _{2}^{2}ds\leq\frac{1}{2}\left\Vert u(t_{0})\right\Vert _{2}^{2}+\int_{t_{0}}^{t}\left(g,u\right)ds.\tag{2.2}
\]
and hence Leray solutions $u$ also satisfy the standard energy inequality.
Such solutions that also satisfy one of the criteria $u\in L^{\theta}((0,T);L^{p}(\Omega),\frac{2}{\theta}+\frac{3}{p}=2,3<p\leq\infty$
are in fact regular solutions of (1.1) on $(0,T)$ by the classic
regularity results mentioned above in the introduction.

We recall that smooth vector fields $v$ which vanish on $\partial\Omega$
in the sense of weak solutions of the Laplace equation satisfy $Av+\nabla(\nabla\bullet v)=\nabla\times(\nabla\times v)$,
and thus if $v$ is divergence-free, i.e. $\nabla\bullet v=0$, then
we have the well-known result that 
\begin{equation}
Av=\nabla\times\nabla\times v\tag{2.3}
\end{equation}
as noted in the introduction. Hence, since under the assumed (e.g.
zero Dirichlet) boundary conditions $A$ is positive definite and
has a well-defined bounded inverse $A^{-1}$ (with $\Vert A^{-1}\Vert_{2}=\lambda_{1}^{-1}$),
 
\begin{equation}
v=A^{-1}(\nabla\times\nabla\times v)\tag{2.4}
\end{equation}
and by applying $A^{s}$ to both sides we obtain Theorem 1 for suitably
smooth $v$. Note that (2.3) holds also in the distributional sense
by considering the application of the appropriate adjoint operators
to smooth test functions; hence (2.4) can hold in this sense for nonsmooth
$v$ as well. In fact $B_{0}\equiv A^{-1}(\nabla\times\nabla\times)$
defines a bounded operator on $H_{0}$. For $D_{c}\equiv\nabla\times$
we have that $B_{0}^{*}=(D_{c}^{*})^{2}A^{-1}$ is clearly a bounded
operator on $H_{0}$ so the result follows by duality; similar arguments
were employed in {[}7{]} to show that the operator $A^{-1/2}Pdiv$
is a bounded operator from $L^{p}(\Omega)$ to $PL^{p}(\Omega)$,
and in {[}1{]} for a related class of operators and spaces.

We begin the proof of Theorem 3 by setting $B_{1}\equiv A^{-1/2}D_{c}$;
by duality again $B_{1}$ is a bounded operator on $L^{p}(\Omega)$,
$1<p<\infty$. Then for $u\in L^{p}(\Omega)$ we have that $A^{s/2}B_{1}u=A^{s/2}(A^{-1/2}D_{c}A^{-s/2})A^{s/2}u=[A^{(s-1)/2}D_{c}A^{-1/2}A^{(1-s)/2}]A^{s/2}u$.
The operator $B_{2}\equiv D_{c}A^{-1/2}$ is clearly a bounded operator
on $L^{p}(\Omega)$, $1<p<\infty$; set $B_{3}\equiv A^{(s-1)/2}D_{c}A^{-1/2}A^{(1-s)/2}=A^{(s-1)/2}B_{2}A^{(1-s)/2}$
then $B_{3}$ is therefore also a bounded operator on $L^{p}(\Omega)$,
$1<p<\infty$, directly if $(1-s)/2\leq0$ and by duality again if
otherwise. Thus $A^{s/2}(B_{1}u)=B_{3}(A^{s/2}u)$ and so $B_{1}\equiv A^{-1/2}D_{c}$
is a bounded operator on $H^{s,p}(\Omega)$ for any real $s$ and
any $p\in(1,\infty)$.

Then from (2.4) we have for any $u\in L^{p}(\Omega)$ that $u=A^{-1}(\nabla\times\nabla\times u)=A^{-1/2}(A^{-1/2}D_{c})(\nabla\times u)=A^{-1/2}B_{1}\omega$.
Since clearly $A^{-1/2}$ is a bounded operator from $H^{s-1,p}(\Omega)$
to $H^{s,p}(\Omega)$ for any real $s$ and any $p\in(1,\infty)$,
we thus have in particular that if $\omega\in H^{s-1,p}(\Omega)$
then $u\in H^{s,p}(\Omega)$ for any real $s$ and any $p\in(1,\infty)$;
setting $s=0$ and $p\in(3,\infty)$ we thus obtain Theorem 3 for
finite $p$ since $H^{0,p}(\Omega)=L^{p}(\Omega)$. For the case $p=\infty$
we observe that since we are on a bounded domain we have for any $r\in[1,\infty)$
that $\Vert A^{-1/2}\omega\Vert_{r}\leq|\Omega|^{1/r}\Vert A^{-1/2}\omega\Vert_{\infty}\leq\Vert A^{-1/2}\omega\Vert_{\infty}$
if $|\Omega|\leq1$ and that $\Vert A^{-1/2}\omega\Vert_{r}\leq|\Omega|^{1/r}\Vert A^{-1/2}\omega\Vert_{\infty}\leq|\Omega|\Vert A^{-1/2}\omega\Vert_{\infty}$
if $|\Omega|\geq1$. Then combining with the remarks above we have
that $\Vert u\Vert_{r}$ is uniformly bounded by $\Vert A^{-1/2}\omega\Vert_{\infty}$
for all $r\in(1,\infty)$, so since $\underset{r\rightarrow\infty}{\lim}\Vert u\Vert_{r}=\Vert u\Vert_{\infty}$
we have that $\Vert u\Vert_{\infty}$ is uniformly bounded by $\Vert A^{-1/2}\omega\Vert_{\infty}$
and we thus obtain Theorem 3 for the case $p=\infty$.

\section{Conclusion}

On reasonable bounded domains with zero boundary conditions $A^{-1}$
is well-defined, and with it we are able to replace the formula $u=-\nabla\times(\nabla^{-1}\omega)$
and the use of the Biot-Savart Law with the identity (2.5). Duality
arguments along the lines of those employed in {[}7{]} and in {[}1{]}
then allow us to use this identity to connect with the standard regularity
criteria $u\in L^{\theta}((0,T);L^{p}(\Omega),\frac{2}{\theta}+\frac{3}{p}=2,3<p\leq\infty$
via suitable operator-theory machinery. The identity (1.3) and the
invertibility of $A$ in fact are the key tools that allow us here
from the outset to consider results of BKM type on bounded domains,
and once in place we see that they allow us to take the extra step
into distributional spaces.


\begin{thebibliography}{10}
\bibitem{key-1} J. Avrin, Global existence and regularity for the
Lagrangian-averaged Navier-Stokes equation with initial data in $H^{1/2}$,
\emph{Comm. Pure Appl. Anal}. \textbf{3} (2004), 353-366.

\bibitem{key-2} J.T. Beale, T. Kato, \& A. Majda, Remarks on the
breakdown of smooth solutions for the 3-D Euler equations, \emph{Commun.
Math. Phys. }\textbf{94} (1984), 61-66.

\bibitem{key-3} L. Escauriaza, G. Seregin, \& V. Sverak, Backward
uniqueness for parabolic equations, \textit{Arch. Rational Mech. Anal.}
\textbf{169} (2003), 147-157.

\bibitem{key-4} J. Fan \& Y. Zhou, A note on regularity criterion
for the 3-D Boussinesq system with partuial viscosity, \textit{Appl}\emph{.
Math. Lett. }\textbf{22} (2009), 802-805.

\bibitem{key-5} J. Fan \& Y. Zhou, A regularity criterion for a fluid
system with the linear Soret effect, \textit{Appl}\emph{. Math. Lett.
}\textbf{25} (2012), 149-152.

\bibitem{key-6}Y. Giga, Solutions for semilinear parabolic equations
in $L^{p}$ and regularity of weak solutions of the Navier-Stokes
system, \emph{J. Diff. Eqs} \textbf{62} (1986), 186\textendash{}212. 

\bibitem{key-7} Y. Giga \& T. Miyakawa, Solutions in $L_{r}$ of
the Navier-Stokes initial-value problem, \textit{Arch. Rational Mech.
Anal.} \textbf{22} (1985), 267-281.

\bibitem{key-8} Z. Grujic \& R. Guberovic, A regularity criterion
for the 3D NSE in a local version of the space of functions of bounded
mean oscillations, \textit{Ann}\emph{. I. H. Poincare,}\textbf{ AN
27} (2010), 773-778.\emph{ }

\bibitem{key-9} H. Kosono \& Y. Taniuchi, Bilinear estimates in BMO
and the Navier-Stokes equations, \emph{Math. Z. }\textbf{235} (2000),
173-194.

\bibitem{key-10} H. Kosono \& Y. Taniuchi, Limiting case of the Sobolev
inequality in BMO, with application to the Euler equations, \emph{Commun.
Math. Phys. Z. }\textbf{214} (2000), 191-200.

\bibitem{key-11} H. Kosono, T Ogawa, \& Y. Taniuchi, The critical
Sobolev inequalities in Besov spaces and regularity criterion to some
semiliniear evolution equations, \emph{Math. Z. }\textbf{242} (2002),
251-278.

\bibitem{key-12} O.A. Ladyzhenskaya, The classical character of generalized
solutions of nonlinear nonstationary Navier-Stokes equations, \textit{Proc.
Steklov Inst. Math}\emph{. }\textbf{92} (1966), 113-131.

\bibitem{key-13} Z. Lei \& Y. Zhou, BKM's criterion and global weak
solutions for magnetohydrodynamics with zero viscosity, \textit{Discrete
Cont. Dyn. Syst}\emph{. }\textbf{25} (2009), 575-583.

\bibitem{key-14} T. Ohyama, Interior regularity of weak solutions
of the time-dependent Navier-Stokes equations, \textit{Proc. Japan
Acad.} \textbf{36} (1960), 273\textendash{}277.

\bibitem{key-15} G. Prodi, Un teorema di unicita per le equazioni
di Navier-Stokes, \emph{Ann. Mat. Pure Appl. }\textbf{48} (1959),
173-182.

\bibitem{key-16} G. Seregin, \& V. Sverak, The Navier-Stokes equations
and backward uniqueness, \textit{Nonlinear problems in mathematical
physics and related topics, II, }\textit{\emph{353\textendash{}366,
Int. Math. Ser. (N. Y.), 2, Kluwer/Plenum, New York, 2002.}}

\bibitem{key-17} W. Rudin, \emph{Functional Analysis,}\textit{\emph{
Kluwer/Plenum, New York, 2002.}}

\bibitem{key-18} J. Serrin, On the interior regularity of weak solutions
of the Navier-Stokes equations, \textit{Arch. Rational Mech. Anal.}
\textbf{9} (1962), 187\textendash{}195.

\bibitem{key-19} W. Von Wahl, Regularity questions for the Navier-Stokes
equations, \textit{Approximation methods for Navier-Stokes problems},
R. Rautmann ed., Lecture Notes in Math. 771, Springer, Berlin-Heidelberg-New
York, 1980, 538-542.

\bibitem{key-20} F.B. Weissler, The Navier-Stokes initial-value problem
in $L^{p}$, \textit{Arch. Rational Mech. Anal.} \textbf{74} (1980),
219-230.

\end{thebibliography}
\end{document}